\documentclass[leqno,11pt]{amsart}
\usepackage{latexsym}
\usepackage{mathrsfs}
\usepackage{amsmath}
\usepackage{amssymb}
\usepackage[bookmarks]{hyperref}
\usepackage{pstricks,pst-plot}

\font\tenscr=rsfs10 % scaled \magstep1
\font\sevenscr=rsfs7 % scaled \magstep1
\font\fivescr=rsfs5 % scaled \magstep1
\skewchar\tenscr='177 \skewchar\sevenscr='177 \skewchar\fivescr='177
\newfam\scrfam \textfont\scrfam=\tenscr \scriptfont\scrfam=\sevenscr
\scriptscriptfont\scrfam=\fivescr
\def\scr{\fam\scrfam}

\font\tenscr=rsfs10 % scaled \magstep1
\font\sevenscr=rsfs7 % scaled \magstep1
\font\fivescr=rsfs5 % scaled \magstep1
\skewchar\tenscr='177 \skewchar\sevenscr='177 \skewchar\fivescr='177
\newfam\scrfam \textfont\scrfam=\tenscr \scriptfont\scrfam=\sevenscr
\scriptscriptfont\scrfam=\fivescr
\def\scr{\fam\scrfam}

\newtheorem{theorem}{Theorem}[section]

\newtheorem{corollary}[theorem]{Corollory}
\newtheorem{lemma}[theorem]{Lemma}
\newtheorem{question}[theorem]{Question}

\newtheorem{prop}[theorem]{Proposition}

\newtheorem{definition} [theorem]{Definition}
\newtheorem{observation} [theorem]{Observation}

\newcommand{\bthm}{\begin{theorem}}
\newcommand{\ethm}{\end{theorem}}
\newcommand{\blem}{\begin{lemma}}
\newcommand{\elem}{\end{lemma}}
\newcommand{\bquest}{\begin{question}}
\newcommand{\equest}{\end{question}}
\newcommand{\bcor}{\begin{corollary}}
\newcommand{\ecor}{\end{corollary}}
\newcommand{\bprop}{\begin{prop}}
\newcommand{\eprop}{\end{prop}}
\newcommand{\bdefn}{\begin{definition}}
\newcommand{\edefn}{\end{definition}}
\newcommand{\bpf}{\begin{proof}}
\newcommand{\epf}{\end{proof}}

\newcommand{\bi}{\begin{itemize}}
\newcommand{\ei}{\end{itemize}}

\newcommand{\bc}{\begin{cases}}
\newcommand{\ec}{\end{cases}}

\newcommand{\ba}{\begin{array}}
\newcommand{\ea}{\end{array}}

\newcommand{\be}{\begin{equation}}
\newcommand{\ee}{\end{equation}}

\newcommand{\bea}{\begin{eqnarray}}
\newcommand{\eea}{\end{eqnarray}}

\newcommand{\beaa}{\begin{eqnarray*}}
\newcommand{\eeaa}{\end{eqnarray*}}

\newcommand{\beastar}{\begin{eqnarray*}}
\newcommand{\eeastar}{\end{eqnarray*}}

\begin{document}

\def\tA{\tilde A}
\def\tX{\tilde X}
\def\tf{\tilde f}
\def\tpi{\tilde \pi}
\def\th{\tilde h}
\def\ta{\tilde \alpha}
\def\th{\tilde h}
\def\tp{\tilde p}
\def\tq{\tilde q}
\def\tl{\tilde l}

\def \vep {\varepsilon}

\def \cd {, \ldots ,}
\def \lf{\| f \|}
\def \bs{\setminus}
\def \ep{\varepsilon}
\def \sig{\Sigma}
\def \si {\sigma}
\def \gam {\gamma}
\def \cinf {C^\infty}
\def \cid {C^\infty (\partial D)}

\def\h#1{\widehat {#1}}
\def\hh#1{\widehat {#1} \bs  {#1}}
\def\hk#1#2{{\widehat {#2}}^{#1}}
\def\hhk#1#2{{\widehat {#2}}^{#1} \bs {#2}}

\def\hr#1{h_r({#1})}
\def\hhr#1{h_r({#1}) \bs {#1}}
\def\hrk#1#2{h_r^{#1}({#2})}
\def\hhrk#1#2{h_r^{#1}({#2}) \bs {#2}}

\def\kh#1#2{{}^{#1}{\widehat {#2}}}
\def\khr#1#2{{}^{#1}h_r({#2})}
\def\kshr#1#2{{}_{#1}h_r({#2})}

\def\<{\langle}
\def\>{\rangle}

\def\spshell{\{ z\in \CN: 1<\|z\| < \rho\}}

\def \C {\mathbb C}
\def \CN {\mathbb C^N}
\def\R {\mathbb R}

\def \Z {\mathbb Z}

\def \bbr {\mathbb R}
\def \bbrs {\mathbb R^2}
\def \bbrn {\mathbb R^n}

\def \cc {\mathscr C}
\def \kk {\mathscr K}
\def \mm {\mathfrak M}

\def \od {\overline D}
\def \ol {\overline L}
\def \oj {\overline J}

\def\sC{{\scr C}}
\def\sF{{\scr F}}
\def\sG{{\scr G}}

%section 3
\def \cn {c_1, \ldots , c_n}
\def \zx { z \in \widehat X }
\def \zxx {\{ z \in \widehat X \} }
\def \siv {\sum^\infty _{j=1}}
\def \sivv {\sum^\infty _{j=n+1}}
\def \smv {\sum^m_{j=n+1}}
\def \snv {\sum^n_{j=j_0+1}}
\def \bjkz {B_{jk} (z_0)}
\def \epp {\varepsilon^2_n/2}
\def \epn {\varepsilon_{n+1}}
\def \nn { {n+1} }
\def \ff {F_{z_{0}}}
\def \fk {f_{k_{1}}}
\def \fkv {f_{k_{v}}}
\def \tf {\tilde{f}}

%---------------------------------------------------
\def \ma {\mathfrak{M}_A}  % or \mathscr *
\def \mb {\mathfrak{M}_B}
\def \mra {\mathfrak{M}_{R(A)}}

\def \xy {(x,y)}
\def \fg {f \otimes g}
\def \vp {\varphi}
\def \pl {\partial L}
\def \ad {A(D)}
\def \aid {A^\infty (D)}
\def \ot {\otimes}
\def \rbl {R_b (L)}
\def\Int {{\rm Int}}

%------------------------ Added by Alex 
\def\rowonly#1#2{#1_1,\ldots,#1_#2}
\def\row#1#2{(#1,\ldots,#1_#2)}
\def\diam{\mathop{\rm diam}\nolimits}

\def\pB{\partial B}
\def\oB{\overline B}

\def\endhat{\widehat{\phantom j}}
\def\endhatk{\widehat{\phantom j}{\phantom |}^k}

\def\pji{p_j^{-1}}
\def\oh{\overline H}

\subjclass[2000]{32E20, 32A65 46J10, 46J15}
\title[Equality of polynomial and rational hulls]{No topological condition implies equality of polynomial and rational hulls}
\author{Alexander J. Izzo}
\address{Department of Mathematics and Statistics, Bowling Green State University, Bowling Green, OH 43403}
\email{aizzo@bgsu.edu}

\begin{abstract}
\vskip 20pt
It is shown that no purely topological condition implies the equality of the polynomial and rational hulls of a set: For any uncountable, compact subset $K$ of a Euclidean space, there exists a set $X$, in some $\C^N$, that is homeomorphic to $K$ and is rationally convex but not polynomially convex.  In addition, it is shown that for the surfaces in $\C^3$ constructed by Izzo and Stout, whose polynomial hulls are nontrivial but contain no analytic discs, the polynomial and rational hulls coincide, thereby answering a question of Gupta.  Equality of polynomial and rational hulls is shown also for $m$-dimensional manifolds ($m\geq 2$) with polynomial hulls containing no analytic discs constructed by Izzo, Samuelsson Kalm, and Wold and by Arosio and Wold.
\end{abstract}

\maketitle

\vskip -2.4 true in
\centerline{\footnotesize\it Dedicated to Hari Bercovici} 
\vskip  2.4 truein

 \section{Introduction}
 
Gabriel Stolzenberg \cite{Stol1} studied conditions for the equality of the polynomial and rational hulls of a set, and raised the question of whether every rationally convex subset of $\C^N$ whose first \v Cech cohomology group with
integer coefficients vanishes is polynomially convex.   He later answered this question in the negative \cite{Stol2} by constructing a set in $\C^2$ that is topologically a disjoint union of two discs and is rationally convex but not polynomially convex.  

In the present paper, we will show that in fact {\it no\/} purely topological condition on a set in $\C^N$ implies the equality of its polynomial and rational hulls.  More precisely, we will show that every uncountable, compact subset of a Euclidean space can be embedded in $\C^N$, for some $N$, so as to be rationally convex but not polynomially convex.  As a complement to this result, we will show that every uncountable, compact subset of a Euclidean space can also be embedded in $\C^N$, for some $N$, so that the polynomial and rational hulls coincide and are nontrivial.  (Here and throughout the paper we say that a hull $X_H$ of a set $X \subset \CN$ is {\it nontrivial} if the set $X_H\setminus X$ is nonempty.)  In both cases, we will show that the embedding can be chosen so that the resulting hulls contain no analytic discs.

It was shown by the present author and Lee Stout \cite[Theorem~1.1]{IzzoS} that every compact surface, with or without boundary, can be embedded in $\C^3$ as a smooth totally real surface having nontrivial polynomial hull containing no analytic discs.  The question of whether the embedded surfaces can be chosen so as to be rationally convex was posed by Purvi Gupta \cite[Section~14]{IzzoS}.  In particular, this raises the question of whether the surfaces constructed by the present author and Stout in the proof of \cite[Theorem~1.1]{IzzoS} are rationally convex.  

Results concerning polynomial convexity often have analogues concerning rational convexity for spaces of dimension one higher.  For instance, compact sets of zero one-dimensional Hausdorff measure are polynomially convex, while compact sets of zero two-dimensional Hausdorff measure are rationally convex.
Therefore, by analogy with the well-known theorem that every nontrivial polynomial hull of a smooth one-dimensional manifold is an analytic variety (see \cite[Corollary~3.1.2]{Stout} for instance), one might expect that every nontrivial rational hull of a smooth two-dimensional manifold must contain an analytic disc.  
On the contrary though, by applying the work of Stolzenberg in \cite{Stol1}, we will show that for the surfaces constructed by the present author and Stout the polynomial and rational hulls in fact coincide.  At the end of \cite{IzzoS}, it is shown that for every compact surface without boundary other than the 
sphere, the real projective plane, and the Klein bottle, an alternative construction yields {\em rationally convex}, totally real embeddings with nontrivial polynomial hulls containing no analytic discs.  Whether the embeddings can be taken to be rationally convex in the case of the sphere, the real projective plane, or the Klein bottle remains open.

Before the work of the present author and Stout discussed above, it was shown by the present author, H\aa kan Samuelsson Kalm, and Erlend 
Forn\ae ss Wold that every manifold, with or without boundary, of dimension $m\geq 3$ can be smoothly embedded in $\C^{2m+4}$ so as to have nontrivial polynomial hull containing no analytic discs \cite{ISW}.  This was recently strengthened by 
Leandro Arosio and Wold \cite{A-W} who showed that the embedding can be obtained in 
$\C^{\lfloor \frac{3m}{2}\rfloor}$.  The proof we will give that the polynomial and rational hulls coincide for the surfaces constructed by the present author and Stout applies also to the manifolds constructed in \cite{ISW} and \cite{A-W}.

The precise statements of the main results of the present paper are as follows.

\begin {theorem} \label{maintheorem1} %Theorem~1.1
Let $K$ be an arbitrary uncountable, compact subspace of $\R^n$. Then there exists a subspace $X$ of $\C^{n+3}$ homeomorphic to $K$ such that $X$ has nontrivial polynomial hull but is rationally convex and satisfies $R(X)=C(X)$.  In $\C^{n+4}$ there exists a subspace $X$ with the above properties such that, in addition, the polynomial hull
$\h X$ of $X$ contains no analytic discs.
\end {theorem}

\bthm  \label{maintheorem2} %Theorem~1.2
Let $K$ be an arbitrary uncountable, compact subspace of $\R^n$. Then there exists a subspace $X$ of $\C^{n+4}$ homeomorphic to $K$ such that the polynomial and rational hulls of $X$ are nontrivial and coincide.  Furthermore, the set $X$ can be chosen so that the hulls contain no analytic discs.
\ethm

\bthm  \label{maintheorem3} %Theorem~1.3
If  $S$ is a smooth compact manifold of dimension $m\geq 2$, with or without boundary, then there exists a $C^\infty$-smooth manifold $\Sigma$ embedded in $\C^{\lfloor \frac{3m}{2}\rfloor}$ that is homeomorphic to $S$ and that has the property that the polynomially convex hull $\h \Sigma$, although  strictly larger than $\Sigma$, contains no analytic discs.  Furthermore, the manifold $\Sigma$ can be chosen to be totally real and so that its polynomial and rational hulls coincide.
\ethm

As discussed above, the manifolds in Theorem~\ref{maintheorem3} were constructed by Arosio and Wold in \cite{A-W}, and in the 2-dimensional case, by the present author and Stout in \cite{IzzoS}.  What is new in Theorem~\ref{maintheorem3} is the final assertion that the polynomial and rational hulls of the constructed manifolds are equal.  Thus our proof of Theorem~\ref{maintheorem3} will be confined to establishing that equality.

Although our focus here is on the relation between the polynomial and rational hulls, in connection with Theorems~1.1 and~1.2, it should be pointed out that the existence, for each uncountable, compact subspace $K$ of $\R^n$, of a homeomorphic copy $X$ in $\C^{n+4}$ such that $\h X$ is nontrivial and contains no analytic discs was proved earlier by the author as \cite[Theorem~1.1]{Izzo}.  This is, of course, contained in each of Theorems~\ref{maintheorem1} and~\ref{maintheorem2} above.  We also point out that a space $K$ is homeomorphic to an uncountable, compact subspace of some Euclidean space if and only if $K$ is a compact metrizable space of finite topological dimension that contains a Cantor set.  This follows immediately from the theorem of Menger and N\"obeling about embedding spaces of finite topological dimension in Euclidean spaces \cite[Theorem~V.2]{HW} and the 
Cantor-Bendixson theorem \cite[Theorem~2A.1]{Mos}. 

One can also consider analogues of Theorems~\ref{maintheorem1} and~\ref{maintheorem2} in an abstract uniform algebra setting.  In fact, Theorems~\ref{maintheorem1} and~\ref{maintheorem2} will be obtained essentially as consequences of such abstract results.  These results will be stated in Section~\ref{definitions} after discussing the relevant notions.

The proofs of Theorems~\ref{maintheorem1} and~\ref{maintheorem2} are based on the existence of certain Cantor sets in complex Euclidean spaces.  In the proof of Theorem~\ref{maintheorem2}, we will use a Cantor set constructed by the author \cite{Izzo} in $\C^3$ whose polynomial and rational hulls are nontrivial and coincide and contain no analytic discs.  The proof of the first half of Theorem~\ref{maintheorem1} will use a Cantor set in $\C^2$ constructed by Vitushkin \cite{Vit} that is rationally convex but whose polynomial hull has interior in $\C^2$.  For the proof of the second part of Theorem~\ref{maintheorem1}, we will need a rationally convex Cantor set with a nontrivial polynomial hull that contains no analytic discs.  We will obtain such a Cantor set from Vitushkin's Cantor set by using a result from the paper \cite{IzzoL} of the author and Norman Levenberg.

\bthm\label{Cantor-rat-convex}
There exists a Cantor set $J$ in $\C^3$ such that $J$ has a nontrivial polynomial hull that contains no analytic discs, $J$ is ration\-ally convex, and $R(J)=C(J)$.
\ethm

In his work on conditions for the equality of polynomial and rational hulls, Stolzenberg introduced a condition that he called the {\em generalized argument principle}.  As discussed in the last section of the present paper, this condition can be thought of as a form of analytic structure in a polynomial hull.  The Cantor set of Theorem~\ref{Cantor-rat-convex} is the first example of a compact set with nontrivial polynomial hull that simultaneously fails to have analytic structure both in the usual sense of the polynomial hull containing no analytic discs and in the sense of the set failing to satisfy the generalized argument principle.

In the next section we briefly recall some standard definitions and notations already used above.  
In Section~\ref{Theorem3} we recall some of Stolzenberg's work on his generalized argument principle and use it to prove Theorem~\ref{maintheorem3}. 
As already mentioned, Section~\ref{definitions} is devoted to abstract uniform algebra analogous of Theorems~\ref{maintheorem1} and~\ref{maintheorem2}.
Theorems~\ref{maintheorem1} and~\ref{maintheorem2} themselves are proved in Section~\ref{polyvsrat}.  In the concluding Section~\ref{conclusion}, we discuss the relationship with the existence of forms of analytic structure in polynomial hulls and raise some open questions.

The inspiration for this paper came from an argument shown to the author by Stout and from the question of Gupta mentioned above.  The author thanks Stout and Gupta for their inspiring correspondence.

It is a pleasure to dedicate this paper to Hari Bercovici who made the year the  author spent as a visitor at Indiana University so enjoyable and stimulating.  Indeed the question posed to the author by Bercovici as to whether a (nonsmooth) one-dimensional manifold can have nontrivial polynomial hull containing no analytic discs (subsequently answered in the affirmative in the author's paper \cite{Izzo}) had a very stimulating effect on the author's research, and in particular, led to the discovery of the Cantor sets with hulls containing no analytic discs used in the present paper.

%%%%%%%%%%%%%%%%%%%%%%%%%%%%%%%%%%%%%%%%%
%                                                    PRELIMINARIES
%
%%%%%%%%%%%%%%%%%%%%%%%%%%%%%%%%%%%%%%%%

\section{Preliminaries}~\label{prelim}
For
$X$ a compact Hausdorff space, we denote by $C(X)$ the algebra of all continuous complex-valued functions on $X$ with the supremum norm
$ \|f\|_{X} = \sup\{ |f(x)| : x \in X \}$.  A \emph{uniform algebra} on $X$ is a closed subalgebra of $C(X)$ that contains the constant functions and separates
the points of $X$.  The maximal ideal space of a uniform algebra $A$ will be denoted by $\ma$.

For a compact set $X$ in $\CN$, the \emph{polynomial hull} $\h X$ of $X$ is defined by
$$\h X=\{z\in\CN:|p(z)|\leq \max_{x\in X}|p(x)|\
\mbox{\rm{for\ all\ polynomials}}\ p\},$$
and the
\emph{rational hull} $\hr X$ of $X$ is defined by
$$\hr X = \{z\in\C^N: p(z)\in p(X)\ 
\mbox{\rm{for\ all\ polynomials}}\ p
\}.$$
An equivalent formulation of the definition of $\hr X$ is that $\hr X$ consists precisely of those points $z\in \C^N$ such that every polynomial that vanishes at $z$ also has a zero on $X$.
The set $X$ is said to be \emph{polynomially convex} if $\h X=X$ and \emph{rationally convex} if $\hr X=X$.

We denote by 
$P(X)$ the uniform closure on $X\subset\CN$ of the polynomials in the complex coordinate functions $z_1,\ldots, z_N$, and we denote by $R(X)$ the uniform closure of the rational functions  holomorphic on (a neighborhood of) $X$. 
Both $P(X)$ and $R(X)$ are uniform algebras, and
it is well known that the maximal ideal space of $P(X)$ can be naturally identified with $\h X$, and the maximal ideal space of $R(X)$ can be naturally identified with $\hr X$.

We denote the open unit disc in the complex plane by $D$.  By an \emph{analytic disc} in $\CN$, we mean an injective holomorphic map $\sigma: D\rightarrow\CN$.
By the statement that a subset $S$ of $\CN$ contains no analytic discs, we mean that there is no analytic disc in $\CN$ whose image is contained in $S$.
An \emph{analytic disc} in the maximal ideal space $\ma$ of a uniform algebra $A$ is, by definition, an injective map $\sigma:D\rightarrow \ma$ such that the function $\hat f\circ \sigma$ is analytic on $D$ for every function $f$ in 
$A$, where $\hat f$ denotes the Gelfand transform of $f$.

%%%%%%%%%%%%%%%%%%%%%%%%%%%%%%%%%%%%%
%%%%%%%%%%%%%%%%%%%%%%%%%%%%%%%%%%%%%

% Proof of Theorem~1.3

%%%%%%%%%%%%%%%%%%%%%%%%%%%%%%%%%%%%%
%%%%%%%%%%%%%%%%%%%%%%%%%%%%%%%%%%%%%

\section{Proof of Theorem~\ref{maintheorem3}}\label{Theorem3}

Except for the final assertion that the manifolds can be chosen so that the polynomial and rational hulls coincide, everything in the statement of Theorem~\ref{maintheorem3} is proven in \cite[Theorem~1.1]{IzzoS} in the case $m=2$, and in \cite[Corollary~1.3]{A-W} for all $m\geq 2$.  Thus we will confine ourselves to proving equality of the polynomial and rational hulls of the  manifolds constructed in \cite{IzzoS} and \cite{A-W} (and noting that the same argument applies also to the manifolds constructed in \cite{ISW}.)  That equality follows readily from the work of Stolzenberg in \cite{Stol1}.   We begin by recalling the relevant results.

The following notion was introduced by Stolzenberg \cite{Stol1}.

\bdefn
A compact set $X$ in $\C^N$ satisfies the {\em generalized argument principle\/} provided that, if $p$ is a polynomial that has a continuous logarithm on $X$, then $0\notin p(\h X)$.
\edefn

Following Stolzenberg will also use the following terminology.

\bdefn
A compact set $X$ is {\em simply coconnected\/} if the first \v Cech cohomology group $\check H^1(X,\Z)$ vanishes.
\edefn

It is well-known that simple coconnectivity of $X$ is equivalent to the
statement that every zero-free function on $X$ has a continuous logarithm.  (See for instance \cite[p.~90]{Gamelin}.)

The following elementary observation is due to Stolzenberg.

\bthm \cite[(2.2)]{Stol1} \label{inclusion}
If $X\subset Y$ are compact sets in $\C^N$ such that $X$ satisfies the generalized argument principle and $Y$ is simply coconnected, then
$\h X \subset \hr Y$.
\ethm

\bcor \cite[Corollary 2.3]{Stol1}
If $X$ is simply coconnected and satisfies the generalized argument principle, then $\h X=\hr X$.
\ecor

The generalized argument principle is of course motivated by the classical argument principle.  The following is a modification of \cite[(2.4)]{Stol1}; since it follows from the classical argument principle in the same manner as \cite[(2.4)]{Stol1}, we omit the proof.

\bprop \label{classicalarg}
If $V$ is an analytic subvariety of an open set of $\C^N$ with $\overline V$ compact, and $p$ is a polynomial with a continuous logarithm on $\partial V=\overline V\setminus V$, then $0\notin p(\overline V)$.
\eprop

Thus the boundary of a relatively compact variety $V$ with $\h {\partial V}=\overline V$ satisfies the generalized argument priniciple.

The following result of Stolzenberg shows that the generalized argument principle is preserved by certain limits.

\bthm \cite[(2.12)]{Stol1}  \label{limit}
Let $(X_n)$ be a sequence of compact sets in $\C^N$ each of which satisfies the generalized argument principle.  If in the Hausdorff metric $X_n\rightarrow X$ and $\h X_n\rightarrow \h X$, then $X$ satisfies the generalized argument principle.
\ethm

We can now proceed with the proof of Theorem~\ref{maintheorem3}

\bpf [Proof of Theorem~\ref{maintheorem3}]
Let $\Sigma$ be one of the surfaces in $\C^3$ constructed in \cite{IzzoS} to prove \cite[Theorem~1.1]{IzzoS} or one of the $m$-manifolds in $\C^{\lfloor \frac{3m}{2}\rfloor}$ constructed in \cite{A-W} to prove \cite[Corollary~1.3]{A-W} (or one of the manifolds constructed in \cite{ISW} to prove \cite[Theorem~1.8]{ISW}).  As noted at the beginning of this section, we need only prove that $\h \Sigma = \hr \Sigma$.
In the case of \cite[Theorem~1.1]{IzzoS}, the surface $\Sigma$ contains a disc $Y$ that contains a compact set $E$ of the form $E=X\times \{0\}$ where $X$ is a set in the distinguished boundary of the bidisc in $\C^2$ constructed by Herbert Alexander in \cite{Alex}.  (See the proof of \cite[Lemma~6.1]{IzzoS}.)  
In the case of \cite[Corollary~1.3]{A-W} or \cite[Theorem~1.8]{ISW}, the construction can be done (and is most naturally done) in such a way that the same statement is true with the disc $Y$ replaced by an $m$-dimensional ball $Y$ and with $E=X\times \{0\}$ replaced by $E=X\times \{0\}^k$ for an appropriate $k$.
Furthermore $\h \Sigma = \h E \cup \Sigma$.  (In the case of \cite[Theorem~1.1]{IzzoS}, see the proof of \cite[Lemma~6.1]{IzzoS} again and \cite[Corollary~3.2]{IzzoS}.)   We will show that $X$ satisfies the generalized argument principle, and hence so does $E$.  Granting this for the moment, Theorem~\ref{inclusion} gives that $\h E$ is contained in $\hr Y\subset \hr \Sigma$.  Consequently, $\h \Sigma = \h E \cup \Sigma\subset \hr \Sigma$, and hence $\h\Sigma=\hr \Sigma$, as desired.

It remains to verify that the set $X$ satisfies the generalized argument principle.  
Alexander \cite{Alex} constructs a sequence $(V_n)$ of varieties in the unit bidisc in $\C^2$.  He notes that $\h{\partial V_n}={\partial V_n} \cup V_n$.  Thus by the observation immediate following Proposition~\ref{classicalarg}, $\partial V_n$ satisfies the generalized argument principle.  One sees from Alexander's proof that in the Hausdorff metric, $\partial V_n\rightarrow X$ and $\h {\partial V_n}\rightarrow \h X$.  Hence, Theorem~\ref{limit} gives that $X$ satisfies the generalized argument principle.
\epf

%%%%%%%%%%%%%%%%%%%%%%%%%%%%%%%%
%%%%%%%%%%%%%%%%%%%%%%%%%%%%%%%%

% Section with the abstract uniform algebra results

%%%%%%%%%%%%%%%%%%%%%%%%%%%%%%%%
%%%%%%%%%%%%%%%%%%%%%%%%%%%%%%%%

\section{Results for abstract uniform algebras} \label{definitions}

In this section we prove analogues of Theorems~\ref{maintheorem1} and~\ref{maintheorem2} in the abstract uniform algebra context.  For that, we first define uniform algebra analogues of the notions of rational hull and rational approximation.  As usual, we denote the Gelfand transform of a function $f$ by~$\hat f$.

\bdefn
Given a uniform algebra $A$ and a compact subset $E$ of the maximal ideal space $\ma$ of $A$, we define the {\em abstract rational hull $\hr E$ of $E$} to be the set
$$\hr E=\{\phi \in \ma: \hat f(\phi)\in \hat f(E)\ {\rm for\ all\ } f\in A \}.$$
We say that $E$ is {\em rationally convex} if $\hr E=E$.
\edefn

An equivalent formulation of the definition of $\hr E$ is that $\hr E$ consists precisely of those points $\phi \in \ma$ such that if $f\in A$ satisfies $\hat f(\phi)=0$, then $\hat f$ has a zero on $E$.  When $E\subset X$ are compact sets in $\C^N$, and $A$ is taken to be $P(X)$, this definition yields the usual rational hull of $E$.

\bdefn
Given a uniform algebra $A$ on a compact space $K$, let $R(A)$ be the uniform closure in $C(K)$ of the set $$\{p/q: p,q\in A\ {\rm and\ } q\ \hbox{\rm is\ 
zero-free}\}.$$  Note that $R(A)$ is a uniform algebra on $K$.
\edefn

When $K$ is a compact set in $\C^N$ and $A=P(K)$, then $R(A)=R(K)$.

As one would expect, there is the following lemma.

\blem \label{mra}
The maximal ideal space of $R(A)$ is the abstract rational hull $\hr K$ of $K$.
\elem

\bpf
A more precise statement of the lemma is that the mapping $\Phi:\mra\rightarrow\ma$ sending each multiplicative linear functional on $R(A)$ to its restriction to $A$ is injective and has range $\hr K$.  For the injectivity of $\Phi$, note that given $\phi\in \mra$ and $g\in A$ zero-free on $K$, we have $\phi(g)\phi(1/g)=\phi(g(1/g))=\phi(1)=1$, so $\phi(1/g)=1/\phi(g)$.  Thus for $p$ and $q$ in $A$ with $g$ zero-free on $K$, we have $\phi(p/q)=\phi(p)/\phi(q)$, so $\phi$ is uniquely determined by its restriction to $A$.  

To see that $\Phi(\mra)$ is contained in $\hr K$, suppose that $\phi$ is in $\Phi(\mra)$ and $f\in A$ satisfies $\hat f(\phi)=0$.  Then $f$ is not invertible in $R(A)$, and hence $f$ has a zero on $K$.  We conclude that $\phi$ is in $\hr K$.

For the reverse inclusion, given $\phi\in \hr K$, we extend $\phi$ to a multiplicative linear functional on the set of quotients of elements of $A$ by setting $\phi(p/q)$ equal to $\phi(p)/\phi(q)$ for $p$ and $q$ in $A$ with $q$ zero-free on $K$.  (We leave to the reader the straight forward verification that this definition does yield a well-defined multiplicative linear functional on quotients. Note that for $q\in A$ zero-free on $K$, we have $\phi(q)\not= 0$ because $\phi$ is in $\hr K$.)  To establish that this extension of $\phi$ extends further to a member of $\mra$, it suffices to verify that $|\phi(p)/\phi(q)|\leq \|p/q\|_K$ for all $p$ and $q$ in $A$ with $q$ zero-free on $K$.  For that, note that for $p$ and $q$ as indicated, $p-[\phi(p)/\phi(q)]q$ is a function in $A$ whose Gelfand transform vanishes at $\phi$ and hence must also vanish at some point $x$ of $K$.  Then $|\phi(p)/\phi(q)|=|p(x)/q(x)|\leq \|p/q\|_K$.
\epf

We are now ready to state and prove uniform algebra analogues of Theorems~\ref{maintheorem1} and~\ref{maintheorem2}.

\bthm\label{abstract1}
Let $K$ be a compact Hausdorff space that contains a Cantor set as a subspace.  Then there exists a uniform algebra $A$ on $K$ such that $\ma$ is strictly larger than $K$ but $K$ is rationally convex and $R(A)=C(K)$.  Furthermore, $A$ can be chosen so that $\ma$ contains no analytic discs.
\ethm

\bthm\label{abstract2}
Let $K$ be a compact Hausdorff space that contains a Cantor set as a subspace.  Then there exists a uniform algebra $A$ on $K$ such that $\ma$ is strictly larger than $K$ and $\hr K=\ma$.  Furthermore, $A$ can be chosen so that $\ma$ contains no analytic discs.
\ethm

It should be remarked that the assertion, contained in each of the above two theorems, that on each compact Hausdorff space $K$ that contains a Cantor set there exists a uniform algebra whose maximal ideal space is strictly larger than $K$ but contains no analytic discs was proved earlier by the author as \cite[Theorem~1.8]{Izzo}.

As discussed in the introduction, the proofs of Theorems~\ref{abstract1} and~\ref{abstract2} are based on the existence of certain Cantor sets in complex Euclidean spaces.  We turn, therefore, to constructing the Cantor set of Theorem~\ref{Cantor-rat-convex}.  For that we need some information about the rational hulls of sets given by a general construction in \cite{IzzoL}.  We recall now the general result of \cite{IzzoL}, establish the need additional information, and then prove Theorem~\ref{Cantor-rat-convex}.

\begin {theorem}\cite[Theorem~1.1]{IzzoL} \label{Wermergen} %Theorem~1.3
Let $X \subset \C^N$ be a compact set whose polynomial hull is nontrivial.  Then there exists a compact set $Y \subset \C^{N+1}$ such that, letting $\pi$ denote the restriction to $\h Y$ of the projection $\C^{N+1} \to \C^N$ onto the first $N$ coordinates, the following conditions hold:
\item {\rm (i)} $\pi (Y) = X$
\item {\rm (ii)} $\pi (\hh Y) = \hh X$
\item {\rm (iii)} $\h Y$ contains no analytic discs
\item {\rm (iv)} each fiber $\pi^{-1} (z)$ for $z \in \h X$ is totally disconnected.
\end {theorem}

\bthm \label{rat-hull-for-Izzo-L}
Let $X$, $Y$, and $\pi$ be as in the preceding theorem.  Then $\pi(\hhr Y)\subset \hhr X$.  In particular, if $X$ is rationally convex, then so is $Y$.   Furthermore, if $R(X)=C(X)$, then $R(Y)=C(Y)$.
\ethm

\bpf
By condition (ii) of the preceding theorem, to show that $\pi(\hhr Y) \subset \hhr X$, it suffices to show that $\pi(\hr Y)\subset \hr X$.  To that end, suppose that $y\in \h Y$ satisfies $\pi(y)\notin \hr X$.  We will show that $y$ is not in $\hr Y$ thereby establishing the desired inclusion.  Because $\pi(y)\notin \hr X$, there is a polynomial $p$ on $\C^N$ such that $p\bigl(\pi(y)\bigr)=0$ but $p$ has no zeros on $X$.  Regard $p$ as a polynomial on $\C^{N+1}$ that is independent of the last variable.  Then $p(y)=0$ and $p$ has no zeros on $\pi^{-1}(X)=Y$.  Thus $y\notin \hr Y$, as desired.

Now suppose that $R(X)=C(X)$.  Then each set of antisymmetry for $R(Y)$ must lie in a single fiber of $\pi$.  Because $Y$ is rationally convex, the maximal ideal space of $R(Y)$ is $Y$, and consequently each maximal set of antisymmetry for $R(Y)$ is connected \cite[p.~119]{St-book}.  Since each fiber of $\pi$ is totally disconnected, we conclude that each maximal set of antisymmetry is a single point.  Consequently, $R(Y)=C(Y)$ by the Bishop antisymmetric decomposition \cite[Theorem~12.1]{St-book}.
\epf

\bpf [Proof of Theorem~\ref{Cantor-rat-convex}]
Take $X$ to be the Cantor set in $\C^2$ whose polynomial hull has interior constructed by Vitushkin in \cite{Vit}, and let $Y$ be the set in $\C^3$ given by then applying Theorem~\ref{Wermergen}.  
Let $J$ be the largest perfect subset of $Y$.  (Recall that a subset of a space is called \emph{perfect} if it is closed and has no isolated points.  Every space contains a unique largest perfect subset (which can be empty), namely the closure of the union of all perfect subsets of the space.)
It follows from conditions (i) and (iv) of Theorem~\ref{Wermergen} that  $J$ is a Cantor set by the well-known characterization of Cantor sets as the compact, totally disconnected, metrizable spaces without isolated points.  
By \cite[Lemma~4.2]{Izzo}, $\h J \bs J \supset \h Y \bs Y$, so condition (ii) of Theorem~\ref{Wermergen} gives that $\h J$ is nontrivial, and since $\h J \subset \h Y$, condition (iii) gives that $\h J$ contains no analytic discs.  
Vitushkin \cite{Vit} noted that $R(X)=C(X)$.  Thus $R(Y)=C(Y)$ by Theorem~\ref{rat-hull-for-Izzo-L}, and consequently $R(J)=C(J)$.  The rational convexity of $J$ is an immediate consequence of the equality $R(J)=C(J)$.
\epf

We are now ready for the proofs of Theorems~\ref{abstract1} and~\ref{abstract2}.

\begin {proof}  [Proof of Theorem~\ref {abstract1}]      
The proof is the same as that of \cite[Theorem~1.8]{Izzo} except that we use the above results to achieve the rational convexity of $K$ and the equality $R(A)=C(K)$.  For completeness, we include the details.

Let $J$ be a Cantor set contained in $K$.   By Theorem~\ref{Cantor-rat-convex}, there exists a uniform algebra $B$ on $J$ such that $\mm_B\bs J$ 
is nonempty but $\mm_B$ contains no analytic discs and $R(B)=C(J)$.  Let $A$ be the uniform algebra on $K$ defined by
$A= \{f\in C(K): f|J\in B\}$.

Each of $K$ and $\mm_B$ can be regarded as subsets of $\mm_A$ in standard ways.  Then $\mm_B$ is the $A$-convex hull of $J$ in $\mm_A$ \cite[Theorem~II.6.1]{Gamelin}, and it follows that $K \cap \mm_B=J$.

Let $\Sigma=K \cup \mm_B$.  We claim that $\Sigma=\mm_A$.  To see this, let $\widehat B$ denote the uniform algebra on $\mm_B$ obtained from $B$ via the Gelfand transform, and let $\widetilde A=\{f\in C(\Sigma):f|\mm_B\in \widehat B\}$.
Then by \cite[Theorem~4]{Bear}, $\Sigma=\mm_{\widetilde A}$.  Since the map $\widetilde A\to A$ given by restriction $(f \mapsto f|K)$ is an isometric isomorphism, this gives that $\Sigma= \mm_A$.

We conclude that if $\mm_A$ contains an analytic disc, then either $\mm_B$ or $\mm_A \bs \mm_B = K \bs J$ must contain an analytic disc.  But $\mm_B$ contains no analytic discs by our choice of $B$, and $K \bs J$ contains no analytic discs because the real-valued functions in $A$ separate points on $K \bs J$.

Now note that $\mm_A \bs K = \mm_B \bs J \ne \emptyset$.

The equality $R(A)=C(K)$ is immediate from the following lemma.  The rational convexity of $K$ is then immediate from Lemma~\ref{mra}.
\end{proof}

\blem
Let $J$ be a Cantor set contained in a compact space $K$, let $B$ be a uniform algebra on $J$, and let $A= \{f\in C(K): f|J\in B\}$.  Then $R(A)= \{f\in C(K): f|J\in R(B)\}$.  In particular, if $R(B)=C(J)$, then $R(A)=C(K)$.
\elem

\bpf
The inclusion $R(A)\subset \{f\in C(K): f|J\in R(B)\}$ is obvious.  For the reverse inclusion, consider an arbitrary function $f\in C(K)$ such that $f|J\in R(B)$, and let $\vep >0$ be arbitrary.  Then there exist functions $p$ and $q$ in $B$ with $q$ zero-free such that
$$\|(p/q)-(f|J)\|<\vep.$$
Set $h=(p/q)-(f|J)$.  By the Tietze extension theorem, $h$ can be extended to a continuous function $\th$ on $K$ such that $\|\th\|_K=\|h\|_J<\vep$.  Set $g=f+\th$.  Then $g$ is in $C(K)$ and $g|J=(f|J)+h=p/q$.  It is not difficult to verify that every zero-free function on a Cantor set has a continuous logarithm.  Let $\alpha$ be such a logarithm for $q$.  Extend $\alpha$ to a continuous function $\tilde\alpha$ on $K$ and set $\tq=e^{\tilde \alpha}$.  Then $\tq$ is an extension of $q$ to a zero-free continuous function on $K$.  Set $\tp=g\tq$.  Then $\tp$ is a continuous extension of $p$ to $K$.  Thus $\tp$ and $\tq$ are in $A$, and 
$$\|f-(\tp/\tq)\|_K=\|f-g\|_K=\|\th\|_K<\vep.$$
Since $\vep>0$ was arbitrary,
we conclude that $f$ is in $R(A)$, as desired.
\epf

Note that the only property of the Cantor set used in the preceding proof was that every zero-free function on the Cantor set has a continuous logarithm, i.e., that the Cantor set is simply coconnected.

\bpf [Proof of Theorem~\ref{abstract2}]
The proof is similar to the proof of Theorem~\ref{abstract1}.  Let $J$ be a Cantor set contained in $K$.  By \cite[Theorem~1.4]{Izzo}, there exists a Cantor set in $\C^3$ with nontrivial polynomial hull such that $P(X)$ contains a dense set of invertible elements.  As noted by Garth Dales and Joel Feinstein \cite{DalesF}, density of the invertible elements in $P(X)$ implies that the polynomial and rational hulls of $X$ coincide and also implies that they contain no analytic discs.  Consequently, there exists a uniform algebra $B$ on $J$ such that $\mb\setminus J$ is nonempty, $\hr J=\mb$, and $\mb$ contains no analytic discs.
As in the proof of Theorem~\ref{abstract1}, define $A$ to be the uniform algebra on $K$ given by
$A= \{f\in C(K): f|J\in B\}$.

Regarding each of $K$ and $\mm_B$ as subsets of $\mm_A$ and arguing exactly as in the proof of Theorem~\ref{abstract1}, we again get that 
%$\ma=K\cup \mb$ and  $K\cap \mb=J$, 
$\ma\setminus K=\mb\setminus J\not=\emptyset$, the maximal ideal space $\ma$ of $A$ contains no analytic discs,  and $\h A|\mb=\h B$.

To see that $\hr K=\ma$, let $x\in \ma\setminus K = \mb\setminus J$ be arbitrary, and let $f\in \h A$ satisfy $f(x)=0$.  Then $f|\mb\in \h B$, so since $\hr J=\mb$ (with respect to $B$), $f|\mb$ must have a zero on $J$.  Thus $f$ has a zero on $K$, so $\hr K=\ma$, as desired.
\epf

%%%%%%%%%%%%%%%%%%%%%%%%%%%%%%%%
%%%%%%%%%%%%%%%%%%%%%%%%%%%%%%%%

% Proofs of Theorems 1.1 and 1.2

%%%%%%%%%%%%%%%%%%%%%%%%%%%%%%%%
%%%%%%%%%%%%%%%%%%%%%%%%%%%%%%%%

\section {Proofs of Theorems~\ref{maintheorem1} and~\ref{maintheorem2}}\label{polyvsrat}

The proofs are similar to the proof of \cite[Theorem~1.1]{Izzo}.
We begin with proving the last part of Theorem~\ref{maintheorem1}.
As in the proof of Theorem~\ref {abstract1}, let $J$ be a Cantor set in $K$, let $B$ be a uniform algebra on $J$ such that $\mb \setminus J$ is nonempty but $\mb$ contains no analytic discs and $R(B)=C(J)$, and let $A=\{ f \in C(K) : f | J \in B \}$. Note that by Theorem~\ref{Cantor-rat-convex}, the uniform
algebra $ B$ can be chosen so as to be generated by three functions $f_1, f_2, f_3$.

Extend each of $f_1, f_2, f_3$ to continuous complex-valued functions $\tf_1, \tf_2, \tf_3$ on $K$.  Let $x_1 \cd x_n$ denote the real coordinate functions 
on $\R^n$.  Choose a continuous real-valued function $\rho$ on $K$ whose zero set is precisely $J$. 
Then the $n+4$ functions $\tf_1, \tf_2, \tf_3, \rho, \rho x_1 \cd \rho x_n$ generate the uniform algebra $A$ by \cite[Lemma~3.8]{ISW} .

Let $F:K\to \C^{n+4}$ be the mapping whose components are the functions  
$\tf_1, \tf_2, \tf_3, \rho, \rho x_1 \cd \rho x_n$, and 
let $X=F(K)$.  Then $P(X)$ is isomophic to $A$ as a uniform algebra and $\mm_{P(X)}$ can be identified with $\h X$.
Since the proof of Theorem~\ref{abstract1} shows that $\ma \setminus K$ is nonempty while $\ma$ contains no analytic discs, and $R(A)=C(K)$, the last part of Theorem~\ref{maintheorem1} follows.  

An alternative approach which avoids the argument in the proof of Theorem~\ref{abstract1} is to apply \cite[Proposition 3.1(i)]{ISW} 
to show directly that, with $G:J\to \C^3$ given by $G(x)=\bigl(f_1(x), f_2(x), f_3(x)\bigr)$, we get
\[ \h X \bs X = \Bigl(\widehat{G(J)} \bs G(J)\Bigr) \times \{0\}^{n+1} \]
Since $G(J)$ is the Cantor set of Corollary~\ref{Cantor-rat-convex}, this also yields the result.  Details are left to the reader.

To see that the example can be found in $\C^{n+3}$ when the requirement that $\h X$ contains no analytic discs is dropped, carry out the above construction starting with Vitushkin's Cantor set \cite{Vit} in $\C^2$  in place of the Cantor set of Theorem~\ref{Cantor-rat-convex}.

The proof of Theorem~\ref{maintheorem2} is essentially the same as the above argument, except that we start with the Cantor set of \cite[Theorem~1.4]{Izzo} in $\C^3$ whose nontrivial polynomial and rational hulls coincide and refer to Theorem~\ref{abstract2} in place of Theorem~\ref{abstract1}.

%%%%%%%%%%%%%%%%%%%%%%%%%%%%%%%%
%%%%%%%%%%%%%%%%%%%%%%%%%%%%%%%%

% Concluding remarks and open questions

%%%%%%%%%%%%%%%%%%%%%%%%%%%%%%%%
%%%%%%%%%%%%%%%%%%%%%%%%%%%%%%%%

\section{Analytic structure and open questions}\label{conclusion}

One sees immediately from Proposition~\ref{classicalarg} that if $X$ is a compact set in $\C^N$ such that every point of $\h X$ lies on a relatively compact analytic variety whose boundary is contained in $X$, then $X$ satisfies the generalized argument principle.  Thus one could think
of the generalized argument principle as a weak form of analytic structure in $\h X$.  Notice, however, that this form of analytic structure\vadjust{\kern 2pt} is incommensurate with the existence of analytic discs in $\h X$.  We have already seen in Alexander's set \cite{Alex} (in the distinguished boundary of the bidisc) a set that satisfies the generalized argument principle but whose polynomial hull contains no analytic discs.  On the other hand Vitushkin's Cantor set fails to satisfy the generalized argument principle but there is an open set of $\C^2$ contained in its polynomial hull.  (Note that for $X$ a Cantor set, the generalized argument principle is equivalent to the statement that $\h X=\hr X$.)  The Cantor set given in Theorem~\ref{Cantor-rat-convex} is the first example of a compact set that simultaneously fails to have an analytic disc in its polynomial hull and fails to satisfy the generalized argument principle.  Each of the compact sets with polynomial hull without analytic discs in \cite{Stol0}, \cite{Wermer1982}, \cite{Alex}, and \cite{DalesF} satisfy the generalized argument principle by the argument given in the second paragraph of the proof of Theorem~\ref{maintheorem3}.  The examples in \cite{Wermerpark}, \cite{DuvalL}, \cite{Izzo0}, and \cite{Izzo} satisfy the generalized argument principle because their polynomial and rational hulls coincide.

The condition that $P(X)$ fails to have a dense set of invertible elements could also be thought of
as a weak form of analytic structure in $\h X$ since this condition holds whenever there is an analytic disc in $\h X$.  Because $X$ obviously satisfies the generalized argument principle whenever $\h X= \hr X$, the elementary observation that $\h X=\hr X$ whenever $P(X)$ has a dense set of invertible elements leads trivially to the following observation about the existence of weak forms of analytic structure.

\begin{observation}\label{trivial}
For every compact set $X\subset\C^N$, there is a weak form of analytic structure in the sense that either $X$ satisfies the generalized argument principle or else $P(X)$ fails to have a dense set of invertible elements.
\end{observation}

Saying that $P(X)$ fails to have a dense set of invertible elements is far from saying that $\h X$ contains analytic structure in the sense of containing an analytic disc.  This observation motivates introducing the following definition and raising the first of the questions below.

\bdefn
A compact set in $\C^N$ satisfies the {\em generalized open mapping principle} if there exists a subset $G$ of $\h X$ such that for every function $f\in P(\h X)$ that is nonconstant on $G$, the set $f(G)$ has interior in $\C$.
\edefn

We can think of the set $G$ in the above definition as a generalized analytic set.

\bquest
Does there exist a compact set $X\subset \C^N$ for which both the generalized argument principle and the generalized open mapping principle fail to hold?  A negative answer would give a more substantive statement about the existence of analytic structure in polynomial hulls than Observation~\ref{trivial}.
\equest

\bquest\label{C2Cantorset}
In the proof of Theorem~\ref{maintheorem2}, we used a Cantor set constructed in \cite{Izzo} lying in $\C^3$ whose polynomial and rational hulls are nontrivial and coincide. Does a Cantor set with these properties exist in $\C^2$?  If so, 
then using that Cantor set in the proof of Theorem~\ref{maintheorem2} would show that the set $X$ of that theorem could be
constructed in $\C^{n+3}$ provided we do not require the hulls of $X$ to contain no analytic discs.  
\equest

A possible candidate for the Cantor set in $\C^2$ requested in Question~\ref{C2Cantorset} is the Cantor set constructed by Walter Rudin in \cite{Rudin-Cantor} (and also presented in \cite[Theorem~III.2.5]{Gamelin} and \cite[pp.~53--54]{Stout}). This is obtained as follows.  Let $E$ be a compact set in the complex plane such that there exists a nonconstant function $f$ continuous on the Riemann sphere $S^2$ and holomorphic off $E$.  (Such a function exists, for instance, whenever $E$ has positive planar measure.)  Choose $f$ so that $f$ has a zero of order one at $\infty$.  Let $g:S^2\rightarrow \C^2$ be given by $g(z)=\bigl(f(z), zf(z)\bigr)$, and let $J=g(E)$.  In the case when $E$ is  a Cantor set, the resulting set $J$ is Rudin's Cantor set.  In all cases, the maximum principle gives immediately that the set $\Sigma=g(S^2)$ is contained in the polynomial hull of $J$.  It can be shown that $\Sigma$ lies also in the rational hull of $J$.  (In case $E$ is simply coconnected (e.g., a Cantor set) this is immediate from Proposition~\ref{classicalarg}.  For the general case see \cite[p.~54]{Stout}.)
Thus the polynomial and rational hulls of $J$ are nontrivial and contain the Riemann surface $\Sigma\setminus J$.
However, neither the rational nor the polynomial hull of this Cantor set has been exactly determined.  Stout (private communication) showed that the set $\Sigma$ has nonzero second \v Cech cohomology group $\check H^2(\Sigma, \C)$ and hence cannot be polynomially convex by a well-known theorem of Andrew Browder \cite{Browder} (or see \cite[p.~93]{Stout}).  Thus the polynomial hull of Rudin's Cantor set must be strictly larger than $\Sigma$, and in fact, by a theorem of Alexander \cite{Alex0} (or see \cite[p.~103]{Stout}) must have topological dimension at least 3 and hence be much larger than $\Sigma$.  

Note that $\Sigma$ is topologically a pinched sphere obtained from the sphere $S^2$ by collapsing the set $f^{-1}(0)$ to a point.
The following simple proof that the cohomology group $\check H^2(\Sigma, \C)$ is nonzero was found by the author. 
Assume without loss of generality that $f^{-1}(0)$ lies outside the open unit disc $D$.  Choose a continuous map $h:S^2\rightarrow S^2$ that is homotopic to the identity and satisfies $h(S^2\setminus D) =\infty$.  For instance, define $h$ by 
$$h(z) =  \begin{cases} \left (\displaystyle{\frac{|z|}{1-|z|}}\right) z & {\rm if \ \ } |z|<1 \cr 
\infty  &  {\rm if \ \ } |z|\geq1 {\rm\ or \ } z=\infty. 
\end{cases} $$
Because $h$ is constant on $f^{-1}(0)$, necessarily $h$ factors through $\Sigma$, that is, there is a continuous map $\tilde h:\Sigma\rightarrow S^2$ such that $h=\tilde h\circ g$.  Then the induced maps on cohomology satisfy $h^*=g^*\circ \tilde h^*$.  Since $h$ is homotopic to the identity, $h^*$ is the identity on $\check H^2(S^2, \C)=\C$.  Therefore, $g^*$ and $\tilde h^*$ are nonzero, and hence, $\check H^2(\Sigma, \C)$ is nonzero.

\bquest
The existence of a  Cantor set in $\C^2$ whose rational hull has interior in $\C^2$ was proven by Henkin in \cite{Henkin}, and the same was shown to hold with $\C^2$ replace by $\C^N$, $N\geq 2$, by the present author in \cite{Izzo}.
 Is there a Cantor set $X$ in some $\C^N$ 
such that $\hr X$ has interior {\em and} $\h X=\hr X$?  Is there an example in $\C^2$?
\equest

\end{document}